\documentclass[12pt,oneside,english]{amsart}
\usepackage[T1]{fontenc}
\usepackage[latin9]{inputenc}
\usepackage{amsthm}
\usepackage{setspace}
\usepackage{amssymb}
\makeatletter
\numberwithin{equation}{section}
\numberwithin{figure}{section}
\theoremstyle{plain}
\newtheorem{thm}{Theorem}

\makeatother

\usepackage{babel}

\begin{document}

\title{\textmd{\textup{\normalsize Chaos in Binary Category Computation}}}

\author{Carlos Pedro Gonçalves}

\maketitle
\begin{center}
Instituto Superior de Ciências Sociais e Políticas, Technical University
of Lisbon 
\par\end{center}

\begin{center}
E-Mail address: cgoncalves@iscsp.utl.pt
\par\end{center}
\begin{abstract}
Category computation theory deals with a web-based systemic processing
that underlies the morphic webs, which constitute the basis of categorial
logical calculus. It is proven that, for these structures, algorithmically
incompressible binary patterns can be morphically compressed, with
respect to the local connectivities, in a binary morphic program.
From the local connectivites, there emerges a global morphic connection
that can be characterized by a low length binary string, leading to
the identification of chaotic categorial dynamics, underlying the
algorithmically random pattern. The work focuses on infinite binary
chains of $\mathcal{C}_{2}$, which is a category that implements
an X-OR-based categorial logical calculus.
\end{abstract}


\section{Introduction}

The first mathematical model of a computer and the framework for what
constitutes mathematical computation theory, introduced by Turing
\cite{key-12}, developed and were informed by what constituted a
logical thinking based upon a paradigm of a classical logic underlying
the rules for discursive validity.

Turing's mathematical computer served a specific purpose, within mathematical
and metamathematical thinking, regarding an algorithmic approach to
mathematical activity. At stake was the ability to build an algorithm
towards an apprehension of a mathematical pattern, the answering of
a mathematical question, the proof of a mathematical theorem, through
mechanical/logical procedures.

\emph{Category computation theory}, on the other hand, deals with
computation based upon categorial structures. The conceptual framework
is different, being closer to what constitutes a different civilizational
moment from that which was lived at the time when Turing proposed
his computer model. In a networked civilization, where the hypertext
and the \emph{world wide web} assumes a fundamental role in the global
interconnectedness, new patterns and new problems can be identified
and worked upon by mathematicians \cite{key-9,key-10}.

While the \emph{classical computer} was a logical machine that had
a classical logical thinking behind it, a thinking that dates back
to Aristotle's logical thinking regarding the formal rules of argument
construction, the \emph{categorial computational structures} are networked
structures with a \emph{hypertextual nature}. It is the discursive
logic of the \emph{hypertext} and of dynamical systemic networks that
underlie the \emph{categorial computational structures}, making category
theory and categorial logic the natural basis for scientific fields
that deal with network situations such as risk science and risk mathematics
\cite{key-7}.

The present work deals with \emph{binary category computation}, dynamical
systems theory and the complexity sciences. The main purpose of the
current work is to address binary patterns produced by \emph{binary
category computational systems} that are morphically compressible
while being (algorithmically) random.

There are two issues that appear interconnected and that allow one
to address this main problem:
\begin{itemize}
\item The emergence of a global morphic connection, with a corresponding
morphic fundament (emergent property), from the local connections;
\item The compressibility of the global morphic connection in what is a
finite string describing how the local connectivities are processed
in the morphic connection between an origin and a target in a morphism.
\end{itemize}
These two issues allow one to address morphic compressibility of what
constitute aleatorial sequences of objects in infinite chains of a
binary category called $\mathcal{C}_{2}$, and comprised of 0 and
1 as objects and of an {}``X-OR'' morphic structure \cite{key-8},
such morphic compressibility, in turn, leads to a result that is obtained
in the probabilistic analysis of such infinite chains: they are deterministic,
following a small length morphic program, and the pattern that emerges
is random-like, such that, if one has only access to the information
on the object sequence in the chain, one is unable to predict the
sequence, beyond a certain probability, that is, the sequence follows
a stochastic process generated by a deterministic program. This constitutes
an example of chaos in \emph{binary category computation}.

All infinite chains in $\mathcal{C}_{2}$, that make emerge an aleatorial
sequence of objects, are chaotic in this sense, which has implications
for number theory, since that any algorithmically random real number,
that is, any real number, whose binary expansion is algorithmically
incompressible, is \emph{morphically compressible}, a result%
\footnote{See \emph{theorem 1}, introduced and proved in \emph{section 3.},
and the discussion that follows.%
} that makes evident the differences between category computation and
classical Turing machine computation, reinforcing, however, a position
sustained from Gödel's work \cite{key-4,key-5}, such that even if
a mathematical pattern has an underlying systemic basis and, thus,
falls within a criterion of (systemic) truth, that truth may be algorithmically
incompressible, that is, inaccessible by finite mechanical means%
\footnote{This statement results from Gödel's \emph{postscriptum} on his and
Turing's works, which appears in the anthology organized by Manuel
Lourenço, \cite{key-5}, a \emph{postscriptum} which was written and
sent by Gödel to be included in the anthology.%
}, but not by what is the systemic computation, which follows the generative
mechanism of the pattern, underlying the fundamental systemic causality
responsible for the pattern itself.

The present work is organized as follows: in \emph{section 2.}, some
elements of category theory are introduced that are relevant to the
present work, as well as the category $\mathcal{C}_{2}$; in \emph{section
3.}, the space of infinite binary strings $2^{\omega}$ is addressed
in its categorial nature and it is connected with chaos in binary
computation and cellular automata theory, through the example of the
Bernoulli map and through the application to the $\mathcal{C}_{2}$
category. It is shown, in \emph{section 3.}, that $\mathcal{C}_{2}$
category's infinite morphic chains can be expressed by a single morphism
in $2^{\omega}$, which allows us to address chaos in such chains,
as well as the main argument of the current work regarding morphic
program size complexity and chaos in $\mathcal{C}_{2}$.

\section{Categories and Binary Category Computation}

\subsection{Categories}

A mathematical category is a mathematical structure composed of objects
and morphisms, the morphisms correspond to binary directional relations
that formalize, within mathematics, a notion of process in which an
\emph{origin }is directionally connected to a \emph{target}, thus,
for two objects, $x$ and $y$, a morphism of fundament $f$ is defined
as:\begin{equation}
x\overset{f}{\rightleftarrows}y\end{equation}
where the top (labeled) arrow corresponds to the \emph{origin arrow},
and connects the origin $x$ to the target $y$, and the bottom arrow
corresponds to the \emph{target arrow}%
\footnote{Although one might use a single-arrow notation, the double-arrow is
preferable when one is addressing the basic systemic motion expressed
by the \emph{morphism}, as argued in\emph{ }\cite{key-6}.%
}, thus, $x$ is directed towards\emph{ $y$}, as \emph{morphic origin},
under the \emph{morphic} \emph{fundament} $f$, and $y$ is directed
towards $x$, as \emph{morphic target}, under the \emph{morphic fundament
$f$}.\emph{ }Morphisms form the building blocks of \emph{morphic
webs}, which are defined as networks of morphisms \cite{key-6}.

A \emph{category $\mathcal{C}$} can be defined as a system of \emph{morphic
webs} such that the \emph{law of morphic identity} is satisfied for
each \emph{morphic node} (\emph{object}), with an associative composition
of morphisms, satisfying a number of relations regarding identity,
as specified next.

Formally, writting $x\in_{ob}\mathcal{C}$, if $x$ is an \emph{object}
of $\mathcal{C}$, the \emph{law of morphic identity} can be expressed
as:\begin{equation}
\forall x\in_{ob}\mathcal{C}:\; x\overset{id_{x}}{\rightleftarrows}x\end{equation}
that is, each object is linked to itself, under its \emph{identity
morphism}. Under composition, for each \emph{morphic chain} of the
sort:\begin{equation}
x\overset{f}{\rightleftarrows}y\overset{g}{\rightleftarrows}z\end{equation}
we have that the composition morphism can be defined as:\begin{equation}
x\overset{g\circ f}{\rightleftarrows}z\end{equation}
If one assumes an algebraic closure for the composition of morphisms%
\footnote{Which is necessary for one to obtain a category.%
}, then, there is a \emph{morphic fundament} $h$, such that:\begin{equation}
h=g\circ f\end{equation}
considering the sign $=$ as an algebraic equality. For composition,
in a category $\mathcal{C}$ , the following algebraic laws are assumed:
\begin{itemize}
\item Identity: \begin{equation}
x\overset{f}{\rightleftarrows}y\Rightarrow f\circ id_{x}=f=id_{y}\circ f\end{equation}

\item Associative:\begin{equation}
x\overset{f}{\rightleftarrows}y\overset{g}{\rightleftarrows}z\overset{h}{\rightleftarrows}w\Rightarrow h\circ\left(g\circ f\right)=\left(h\circ g\right)\circ f\end{equation}

\end{itemize}
With these algebraic laws, a \emph{mathematical category} can be addressed
as a system of \emph{computational webs}, whose building blocks (elementary
computations) are the morphisms.

\subsection{Behavior of the automorphism fundament regarding the identity morphism
- emergence and automorphisms}

Let us, first, recover the law of identity over algebraic composition:\begin{equation}
x\overset{f}{\rightleftarrows}y\Rightarrow f\circ id_{x}=f=id_{y}\circ f\end{equation}
This law is general for two different objects, however, one may consider
two types of morphic structures that satisfy the categorial axioms
but that behave differently in regards to this law. The first one
is:\begin{equation}
\forall x\in_{ob}\mathcal{C}:\; x\overset{f}{\rightleftarrows}x\Rightarrow f\circ id_{x}=f=id_{x}\circ f\end{equation}
The systemic interpretation of this rule is that the automorphism
is not sistemically contained in the identity morphism. In the morphic
connection of the object to itself, the property of the connection
$f$ is an emergent property that is not implicated by the identity
morphism. In the Finset we have a trivial example of this type of
morphic structures.

Now, in the second type of morphic structure, the identity connection
already brings with it all of the automorphisms, that is, we have:\begin{equation}
\forall x\in_{ob}\mathcal{C}:\; x\overset{f}{\rightleftarrows}x\Rightarrow f\circ id_{x}=id_{x}=id_{x}\circ f\end{equation}
In such a category, the law of identity over algebraic composition
holds for any morphism between different objects, but not for the
automorphisms. One may also have a combination of the two types of
morphic structures, such that, for some automorphism, (2.9) holds,
while, for others, (2.10) holds.

The second type of morphic structures takes place in $\mathcal{C}_{2}$,
which is now introduced.

\subsection{Binary category computation and $\mathcal{C}_{2}$}

The simplest computational structure, that of a \emph{single bit}
binary computing web, can be built by considering a binary alphabet
as an object base for a binary category $\mathcal{C}_{2}$, comprised
by the numbers $0$ and $1$ as objects, and by the following non-identitarian
singular morphisms:\begin{equation}
\forall x\in_{ob}\mathcal{C}_{2}:\; x\overset{\hat{0}}{\rightleftarrows}x\end{equation}
\begin{equation}
\forall x\in_{ob}\mathcal{C}_{2}:\; x\overset{\hat{1}}{\rightleftarrows}1-x\end{equation}
with the operator algebraic structure, with respect to composition,
being defined in terms of the operator product such that, for $\hat{u}$
ranging in the morphic fundament domain comprised of $\hat{0}$, $\hat{1}$,
$id_{x}$, with $x=0,1$, we have:\begin{equation}
\hat{0}\circ\hat{u}=\hat{u}\circ\hat{0}=\hat{u}\end{equation}
\begin{equation}
\hat{u}\circ\hat{1}=\hat{1}\circ\hat{u}=\hat{1},\;\hat{u}\neq\hat{1}\end{equation}
\begin{equation}
\hat{0}^{2}=\hat{0}\circ\hat{0}=\hat{0}\end{equation}
\begin{equation}
\hat{1}^{2}=\hat{1}\circ\hat{1}=\hat{0}\end{equation}
With these relations, the operator algebraic structure satisfies the
identity laws as well as the associative laws \cite{key-8}. However,
regarding the \emph{automorphism} $x\overset{\hat{0}}{\rightleftarrows}x$,
$\mathcal{C}_{2}$ satisfies (2.10) rather than (2.9), hence corresponding
to a category of the second type, where the \emph{identity morphism}
implicates the \emph{automorphism}.

Now, given the elementary morphisms in $\mathcal{C}_{2}$, it follows
that $\mathcal{C}_{2}$ can be addressed as a binary computational
system of \emph{morphic webs}, whose objects are \emph{bits} and the
non-identity \emph{morphisms} correspond, respectively, to an operation
that leaves the \emph{bit} unchanged and a Boolean negation.

If one considers the input states $\left(\hat{u}x\right)$, where
$x$ is the \emph{origin}, we have the following Boolean table for
the target:

\[
\begin{array}{cccc}
(\hat{0}0) & (\hat{0}1) & (\hat{1}0) & (\hat{1}1)\\
0 & 1 & 1 & 0\end{array}\]
which is an X-OR Boolean function. Thus, $\mathcal{C}_{2}$ is called
an \emph{X-OR-categorial computational structure}. Besides $\mathcal{C}_{2}$,
there are fifteen other categories with the same objects as $\mathcal{C}_{2}$
and with two non-identitarian morphisms, but whose two morphisms lead
to each of the fifteen other Boolean functions. In the present work,
of the sixteen alternative structures, we focus only on $\mathcal{C}_{2}$.

\subsection{Infinite morphic chains in $\mathcal{C}_{2}$}

Let $\mathfrak{M}$ be an infinite morphic chain in $\mathcal{C}_{2}$,
with the general form:\begin{equation}
s_{0}\overset{\hat{u}_{1}}{\rightleftarrows}s_{1}\overset{\hat{u}_{2}}{\rightleftarrows}s_{2}...\end{equation}
the chain is completely determined by the initial condition $s_{0}$
and an infinite sequence of operators $\hat{u}_{1},\hat{u}_{2},...$
. Since the operators can assume only binary values ($\hat{0}$ or
$\hat{1}$), it follows that one may concatenate the infinite sequence
$\hat{u}_{1},\hat{u}_{2},...$ in one single binary string $\mathfrak{P}$
characterizing what we call a morphic program for the chain:\begin{equation}
\mathfrak{P}\equiv\hat{u}_{1}\hat{u}_{2}...\end{equation}
likewise, one may concatenate the objects in a single string, expressing
the object sequence:\begin{equation}
S\equiv s_{0}s_{1}s_{2}...\end{equation}

Addressing, now, the object sequence $S$, there are two possible
configurations of infinite chains in $\mathcal{C}_{2}$ that can be
considered:
\begin{enumerate}
\item \emph{Eventually periodic: }The chain tends to a repeating finite
cycle of objects which is described by a finite string $w_{cycle}$
in the binary alphabet;
\item \emph{Nonperiodic: }The chain never tends to a repeating finite cycle,
being nonperiodic.
\end{enumerate}
In order to define eventual periodicity we need to address a distance
notion, in regards to infinite strings. To that end, let {}``$\sqsubset$''
denote the prefix relation, such that, for any finite binary string
$w$, of length $|w|$, $w\sqsubset X$, where $X$ is some infinite
binary string, if, and only if, the first $|w|$ bits of $X$ coincide
with $w$. With this relation, the following metric is introduced:\begin{equation}
d\left(X,Y\right)\equiv2^{-n}\end{equation}
where $X$, $Y$ are two infinite strings and $n$ is such that:\begin{equation}
n=\max\left\{ \left|w\right|:\: w\sqsubset X\wedge w\sqsubset Y\right\} \end{equation}
if the two strings coincide $d\left(X,Y\right)=2^{-\infty}=0$. Thus,
$n$ is the size of the largest common prefix between $X$ and $Y$.

Now, given the object sequence of an infinite chain $S$ in $\mathcal{C}_{2}$,
let us define the following sequence:\[
S_{0}\equiv S=s_{0}s_{1}s_{2}...,\: S_{1}=s_{1}s_{2}...,\:...,\: S_{n}=s_{n}s_{n+1}...\]

Let now, $S_{cycle}\left[w\right]$ result from the infinite concatenation
of a finite string $w$ with itself, thus, corresponding to an infinite
repetition of the string $w$. The chain in $\mathcal{C}_{2}$, with
object sequence $S$, is, then, stated to be eventually periodic,
if, and only if, there is a $S_{cycle}\left[w\right]$ such that $d\left(S_{n},S_{cycle}\left[w\right]\right)\rightarrow0$,
with $n\rightarrow+\infty$. For a chain $\mathfrak{M}$ that is eventually
periodic with cycle $w$, $\left|w\right|$ is called its cycle number
and we write:\begin{equation}
CN\left(\mathfrak{M}\right)=\left|w\right|\end{equation}

Now, for the second case, of a chain that is not eventually periodic,
their cycle number must be infinite and coincident with Cantor's first
transfinite number:\begin{equation}
CN\left(\mathfrak{M}\right)=\aleph_{0}\end{equation}

It is for chains that are not eventually periodic that one may find
chaotic dynamics in $\mathcal{C}_{2}$, to address such a dynamics,
and more general binary dynamical structures in categories, it becomes
useful to introduce a class of \emph{binary categorial automata} in
the space of infinite binary strings $2^{\omega}$.

\section{Binary Categorial Automata in the Space $2^{\omega}$, Morphic Program
Size Complexity and Chaos}

\subsection{Binary categorial automata in $2^{\omega}$}

Let $2^{\omega}$ be the space of infinite binary strings, and address
it as a category by considering its (morphic) points as objects and
the maps between points as morphisms. A map between two points $S$,
$S'$ being defined as follows:\begin{equation}
S\overset{f}{\rightleftarrows}S'\end{equation}
such that there is a local string-level element morphism: \begin{equation}
s_{n}\overset{f_{n}}{\rightleftarrows}s'_{n}\end{equation}
for each $n\geq0$. Thus, $S\overset{f}{\rightleftarrows}S'$ can
be considered as a composite morphism with fundament $f$, whose connection
takes place at each bit, in each string, through a local morphic connection,
with fundament $f_{n}$. In systemic terms, and considering the binary
relational nature of the morphism, the morphic fundament $f$ can
be considered as an emergent property resulting from the local morphic
connections. Computationally, one may consider to be dealing with
an automaton that, given the input $S$, transforms it into the output
$S'$, morphically linking the two morphic points of $2^{\omega}$.

The automaton framework can be taken further, in such a way that there
is a familiarity with cellular automata, in particular, this computational
system is equivalent to a cellular automaton system with a right-infinite
array of cells. Indeed, local configurations of rules can be considered
that replicate cellular automaton rules. As an example, one may consider
the Rule 170 which, applied to a right-infinite array of cells, implements
the Bernoulli map \cite{key-15}:\begin{equation}
x_{n}=B\left(x_{n-1}\right)=2x_{n-1}\: mod\:1\end{equation}
This map is implemented in $2^{\omega}$ as follows:\begin{equation}
S\overset{B}{\rightleftarrows}S'\end{equation}
with the local morphisms:\begin{equation}
s_{n}\overset{s_{n+1}}{\rightleftarrows}s_{n}'=s_{n+1}\end{equation}
such that we have the following morphic structure for a single iteration:\begin{equation}
\begin{array}{cccc}
s_{0} & s_{1} & s_{2} & ...\\
_{s_{1}}\downarrow\uparrow & _{s_{2}}\downarrow\uparrow & _{s_{2}}\downarrow\uparrow\\
s_{1} & s_{2} & s_{3} & ...\end{array}\end{equation}
effectively implementing the Bernoulli shift map on $S$. A simplifying
approach was taken here, in defining the local morphisms, that demands
some explanation. In each local connection, the fundament provides
the information needed to assign to each origin bit $s_{n}$ a target
$s_{n}'$ which, in this case, is $s_{n+1}$, the $n+1$-th element
of $S$. Thus, each local morphic connection has the form $s_{n}\overset{s_{n+1}}{\rightleftarrows}s{}_{n}'=s_{n+1}$.
These local connections can be used to produce a string composed of
substrings of the form $s_{n}s_{n+1}s_{n}'$ that synthesize the basic
transformation in terms of \emph{input/fundament/output}, this leads
to a small length string describing the automaton's program, and given
by: \begin{equation}
000011100111\end{equation}
the string should be read from left to right and three bits at the
time, in accordance with $s_{n}s_{n+1}s_{n}'.$ In terms of Boolean
function tables, (3.6) corresponds to a \emph{Boolean transfer function}
which transfers, or copies, the value of $s_{n+1}$to the output $s_{n}'$
and ignores the value of $s_{n}$.

Another alternative would have been to use the cellular automaton
notational rules, but we use the strings above in order to be able
to address computational complexity issues related to the size of
binary strings describing morphic programs, and the relation between
these programs and the local connectivities.

To implement the iterations of the Bernoulli shift map, the morphic
program, described by the string (3.7), is applied iteratively, so
that we have morphic chains, in $2^{\omega}$, for the iterations
of the Bernoulli map, expressed in terms of:\begin{equation}
S_{0}\overset{B}{\rightleftarrows}S_{1}\overset{B}{\rightleftarrows}S_{2}...\end{equation}
While, for the Bernoulli map, the iterations are obtained through
morphic chains in $2^{\omega}$, one may recover, within $2^{\omega}$,
through a single morphism, the systemic process underlying the infinite
chains in $\mathcal{C}_{2}$. An infinite chain in $\mathcal{C}_{2}$
can be obtained as a single morphism in $2^{\omega}$ by, first, rewritting
the operator string $\mathfrak{P}$, defined in (2.17), as a binary
string without the {}``hats'', and, therefore, representing it as
a point in $2^{\omega}$, and, second, by defining the morphism in
$2^{\omega}$ from this operator string to the object sequence string
as: \begin{equation}
\mathfrak{P}\overset{\pi}{\rightleftarrows}S\end{equation}
where $\pi$ has the following general properties, for the local morphic
connections, at each string point:\begin{equation}
u_{n}=0\overset{s_{n-1}=s}{\rightleftarrows}s_{n}=s\end{equation}
\begin{equation}
u_{n}=1\overset{s_{n-1}=s}{\rightleftarrows}s_{n}=1-s\end{equation}
with $s=0,1$. For $n=0$, the morphic connection is given by:\begin{equation}
u_{0}=0\overset{s_{0}}{\rightleftarrows}s_{0}\end{equation}
In this way, we have added a $0$ to our operator string which, thus,
is given by: $u_{0}u_{1}u_{2}...$, leaving, once more, the {}``hats''
behind.

With this morphism, each element of $S$, $s_{n}$ (with $n>0$) is
obtained through the action of $u_{n}$ over $s_{n-1}$, an action
that is in accordance with the operator rules for $\mathcal{C}_{2}$,
via the morphic connection with fundament $\pi,$ which involves the
local connections (3.10) and (3.11).

Now, following the same approach that was introduced for the Bernoulli
shift map (\emph{input/fundament/output}) the binary string that expresses
the automaton's computation is given by:\begin{equation}
000011101110\end{equation}
this string synthesizes all of the local morphic connections.

Thus, as in the case of the Bernoulli maps, we were able to provide
for a small length rule, namely, in both cases, it was obtained a
finite morphic program string of length twelve, which expresses the
underlying rule that transforms the input string in the output string.

The above string can, in turn, be recast within $\mathcal{C}_{2}$,
in the form of a finite morphic chain:\begin{equation}
0\overset{\hat{0}}{\rightleftarrows}0\overset{\hat{0}}{\rightleftarrows}0\overset{\hat{0}}{\rightleftarrows}0\overset{\hat{1}}{\rightleftarrows}1\overset{\hat{0}}{\rightleftarrows}1\overset{\hat{0}}{\rightleftarrows}1\overset{\hat{1}}{\rightleftarrows}0\overset{\hat{1}}{\rightleftarrows}1\overset{\hat{0}}{\rightleftarrows}1\overset{\hat{0}}{\rightleftarrows}1\overset{\hat{1}}{\rightleftarrows}0\end{equation}

There is an akinness between the result for the Bernoulli shift map
and the result for $\mathcal{C}_{2}$'s infinite length chains, in
that the later fall under the same systemic typology as systems capable
of, given certain inputs, to produce a chaotic dynamics.

In the case of $\mathcal{C}_{2}$, this dynamics is present in the
morphic chain itself, and can be seen to relate to the distribution
of the morphic fundaments, that is, to $\mathfrak{P}$, while, in
the case of the Bernoulli shift map, it is present in the structure
of the initial state $S_{0}$. In the $2^{\omega}$ categorial automata
representation both cases of chaos are linked to the interplay between
the structure of the initial condition%
\footnote{$\mathfrak{P}$, for $\mathcal{C}_{2}$, and $S_{0}$, for the Bernoulli
map.%
} and the local morphic connection that make emerge a global morphic
connection.

Thus, computationally, the morphism falls under Varela's notion of
\emph{enactive computation} \cite{key-13,key-14}, it is this \emph{enactive
nature} of morphic\emph{ }computation that leads to a global morphic
compressibility in a small length morphic program, even when one is
dealing with algorithmically random structures, providing for an expansion
of the notion of chaos to encompass those examples that fall within
the algorithmic incompressible typology.

The above result for $\mathcal{C}_{2}$ shows that even if the dynamics
produces a pattern that, considered independently of the generative
systemic fundament, is algorithmically incompressible, once the generative
systemic fundament is addressed, one recovers, within a mathematical
setting, a morphic computational compressibility, expressing a small-length
rule that explains how an input generates an output.

While, in the Bernoulli shift map, this rule acts locally at the bits
of each state, in $\mathcal{C}_{2}$, the rule acts locally at each
morphic connection in the chain, but both rules are equally compressible,
providing a picture of low-dimensional chaos, within a binary computational
framework, in connection with a systemic generated randomness, whose
order is accessible and expressible in a small (morphic) program,
twelve bits long.

The morphic compressibility of algorithmically incompressible infinite
binary strings follows, as a corollary, from the following theorem:
\begin{thm}
All infinite binary strings are morphically compressible, under the
categorial automaton \textup{000011101110}.\end{thm}
\begin{proof}
Let $\mathcal{C}_{2}^{\infty}$ be a category whose objects are the
infinite chains of $\mathcal{C}_{2}$. Let a morphism in $\mathcal{C}_{2}^{\infty}$
be defined as follows: given $C$ and $C'\in_{ob}\mathcal{C}_{2}^{\infty}$,
$C\overset{f}{\rightleftarrows}C'$, with the \emph{fundament} $f$
defined such that, for $n\geq0$, the point $s_{n}$ of the origin
chain $C$ is mapped to the corresponding point $s'_{n}$ in the target
chain $C'$, under the local morphism:\begin{equation}
s_{n}\overset{f_{n}}{\rightleftarrows}s'_{n}\end{equation}
Now, for $\mathcal{C}_{2}^{\infty}$ with such a structure, one may
define the functor:\begin{equation}
\mathcal{C}_{2}^{\infty}\overset{\varphi}{\rightleftarrows}2^{\omega}\end{equation}
such that, for each $C\in_{ob}\mathcal{C}_{2}^{\infty}$, $\varphi\left(C\right)=S$,
where $S$ is the object sequence of the chain $C$, $\varphi\left(id_{C}\right)=id_{\varphi(C)}$,
and $\varphi\left(f\right)=f$, such that $\varphi\left(C\right)\overset{f}{\rightleftarrows}\varphi\left(C'\right)$,
where the morphic fundament $f$ results, in $2^{\omega}$, from the
morphism with local morphic connections defined by (3.15) at each
string point, finally, regarding composition, $\varphi\left(g\circ f\right)=\varphi\left(g\right)\circ\varphi\left(f\right)=g\circ f$.

Thus, the functor maps the infinite chains to their corresponding
object sequences and each morphism in $\mathcal{C}_{2}^{\infty}$
to the corresponding morphism in $2^{\omega}$, defined by the same
local morphic connections. Composition and identity are mapped such
that the local morphic structure defined by (3.15) is preserved.

Now, given any $C$ in $\mathcal{C}_{2}^{\infty}$, with object sequence
$S$, it is the only chain with such object sequence, which comes
from the fact that there is, in any chain of $\mathcal{C}_{2}$ (finite
or infinite), one, and only one, possible morphic fundament sequence
for each object sequence, a property which results from the fact that,
in $\mathcal{C}_{2}$, the origin and target completely determine
the fundament of non-identitarian morphisms. Since any $C\in_{ob}\mathcal{C}_{2}^{\infty}$,
with object sequence $S$, is the only object of $\mathcal{C}_{2}^{\infty}$
with such object sequence, we have:\begin{equation}
\varphi\left(C\right)=\varphi\left(C'\right)=S\Leftrightarrow C=C'\end{equation}
which means that the functor is one-to-one, in the objects. It is
also one-to-one in the morphisms, since the morphisms are defined
from their local connectivities and the corresponding \emph{origin}
and \emph{target}. On the other hand, given any $S\in_{ob}2^{\omega}$,
there is one, and only one, $C\in_{ob}\mathcal{C}_{2}^{\infty}$,
such that $\varphi\left(C\right)=S$, a property that results from
the fact that, given any infinite string, it can be converted into
one, and only one, infinite chain in $\mathcal{C}_{2}$, thus, all
of the objects of $2^{\omega}$ have a corresponding origin in $\mathcal{C}_{2}$,
under the functor with fundament $\varphi$, which means that the
functor is onto in the objects. It is also onto in the morphisms,
since, for any morphism $S\overset{f}{\rightleftarrows}S'$, there
is one, and only one, morphism in $\mathcal{C}_{2}^{\infty}$, $C\overset{f}{\rightleftarrows}C'$
such that $\varphi\left(C\right)=S\overset{f}{\rightleftarrows}\varphi\left(C'\right)=S'$,
which means that all of the morphisms of $2^{\omega}$ have a corresponding
origin in $\mathcal{C}_{2}$.

Thus, $\mathcal{C}_{2}^{\infty}$ and $2^{\omega}$ are isomorphic,
under $\varphi$, with inverse functor $\varphi^{-1}$. This means
that, given any infinite binary string $S\in_{ob}2^{\omega}$, one
can pass from the string representation to the morphic chain representation,
through $\varphi^{-1}\left(S\right)$ and, then, given the morphic
program $\mathfrak{P}_{\varphi^{-1}(S)}$, underlying the infinite
chain $\varphi^{-1}(S)$, introduce the morphism in $2^{\omega}$:\begin{equation}
\mathfrak{P}_{\varphi^{-1}(S)}\overset{\pi}{\rightleftarrows}S\end{equation}
which corresponds to the action of the categorial automaton 000011101110
on the input string $\mathfrak{P}_{\varphi^{-1}(S)}$, leading to
the output string $S$. Thus, the generative mechanism of $S$ is
completely determined by the morphism (3.18) with morphic program
string 000011101110, which proves that any infinite binary string
can be morphically compressed.
\end{proof}
Morphic compressibility does not coincide with Turing-based algorithmic
compressibility, because the categorial automaton 000011101110 is
able to work with infinite strings as inputs and produce an output,
following the local connectivity rules in a non-sequencial way, that
is, each local morphic connection occurs simultaneously in the morphism
$\mathfrak{P}_{\varphi^{-1}(S)}\overset{\pi}{\rightleftarrows}S$.

Since any infinite string can be morphically compressed, even an algorithmically
incompressible infinite binary string is morphically compressible,
under the morphism (3.18), a result that immediately follows from
the above theorem, allowing us to address algorithmically random and
random-like infinite binary strings%
\footnote{Since not all random-like infinite binary strings are algorithmically
random \cite{key-2}, we refer the two cases explicitly, allowing
for different criteria for randomness.%
} as chaotic. This matter is, now, further explored.

\subsection{Chaos in $\mathcal{C}_{2}$}

Let us consider the categorial automaton 000011101110 in $2^{\omega}$,
for $\mathcal{C}_{2}$ chains. The first matter to be relevated is
that each element $s_{n}$ of the output string is connected to the
previous element $s_{n-1}$ and to the input morphic fundament operator
string element $u_{n}$, such that, the binary code for the automaton
is comprised of the formal language:

\begin{equation}
\mathfrak{L}_{\mathcal{C}_{2}}=\left\{ 000,011,101,110\right\} \end{equation}
where each string is in the form: \begin{equation}
u_{n}s_{n-1}s_{n}\end{equation}

Based on the binary code, one can introduce a probabilistic analysis
of the output string, given by the conditional probability measures,
for each string element $s_{n}$, conditional on an available information
structure $I_{n}$:

\begin{equation}
\mathbb{P}\left[s_{n}=s|I_{n}\right],\:\mathbb{P}\left[s_{n}=1-s|I_{n}\right]\end{equation}
If $I_{n}$ includes both $s_{n-1}$ and $u_{n}$ one is able to know
the $n$-th string element, thus, for the $n$-th object of the infinite
chain in $\mathcal{C}_{2}$, the conditional probability measures
are as follows:

\begin{equation}
\begin{cases}
\mathbb{P}\left[s_{n}=s|s_{n-1}=s,u_{n}=0\right]=1\\
\mathbb{P}\left[s_{n}=1-s|s_{n-1}=1-s,u_{n}=0\right]=0\end{cases}\end{equation}
and\begin{equation}
\begin{cases}
\mathbb{P}\left[s_{n}=1-s|s_{n-1}=s,u_{n}=1\right]=1\\
\mathbb{P}\left[s_{n}=s|s_{n-1}=1-s,u_{n}=1\right]=0\end{cases}\end{equation}

This result is in keeping with the fact that we are dealing with a
deterministic automaton, where each step is determined in accordance
with the program string 000011101110. Now, let us assume that the
information of the $\mathfrak{P}$ string is not given to the probabilistic
analysis, such that the only available information that is given is
the object string itself. Then, the probabilities depend upon the
statistical structure of the string.

Let us assume, for instance, that the output string $S$, that provides
for the object sequence of the morphic chain in $\mathcal{C}_{2}$,
corresponds to one of the points in $2^{\omega}$ that is statistically
equivalent to an infinite string of independent Bernoulli (random)
trials, such that there is a constant probability (in the frequencist
sense) of $s$ occurring at any place in the string. Then, we have,
through statistical independence:

\begin{equation}
\begin{cases}
\mathbb{P}\left[s_{n}=s|s_{n-1}\right]=\mathbb{P}\left[s_{n}=s\right]=P_{s}\\
\mathbb{P}\left[s_{n}=1-s|s_{n-1}\right]=\mathbb{P}\left[s_{n}=1-s\right]=1-P_{s}\end{cases}\end{equation}

This only takes place under certain conditions for $\mathfrak{P}$,
such that we have the following probabilities:\begin{equation}
\mathbb{P}\left[u_{n}=0|s_{n-1}=s\right]=P_{s}=\mathbb{P}\left[u_{n}=1|s_{n-1}=1-s\right]\end{equation}
\begin{equation}
\mathbb{P}\left[u_{n}=0|s_{n-1}=1-s\right]=1-P_{s}=\mathbb{P}\left[u_{n}=1|s_{n-1}=s\right]\end{equation}
Thus, there is a statistical dependence between the $u_{n}$ and the
$s_{n-1}$, such that, whenever, $s_{n-1}=s$, there is a probability
of $P_{s}$ that $u_{n}=0$, leaving the bit unchanged ($s_{n}=s_{n-1}=s$),
on the other hand, whenever $s_{n-1}=1-s$, then, there is also a
probability of $P_{s}$ that $u_{n}=1$, changing the bit from $s_{n-1}=1-s$
to $s_{n}=s$. A similar reasoning applies to $s_{n-1}=1-s$. Putting
it all together the following three transition matrices are obtained:\begin{equation}
\left(\begin{array}{cc}
\mathbb{P}\left[s_{n}=0|s_{n-1}=0\right]=P_{0} & \mathbb{P}\left[s_{n}=1|s_{n-1}=0\right]=1-P_{0}\\
\mathbb{P}\left[s_{n}=0|s_{n-1}=1\right]=P_{0} & \mathbb{P}\left[s_{n}=1|s_{n-1}=1\right]=1-P_{0}\end{array}\right)\end{equation}
\begin{equation}
\left(\begin{array}{cc}
\mathbb{P}\left[u_{n}=0|s_{n-1}=0\right]=P_{0} & \mathbb{P}\left[u_{n}=1|s_{n-1}=0\right]=1-P_{0}\\
\mathbb{P}\left[u_{n}=0|s_{n-1}=1\right]=1-P_{0} & \mathbb{P}\left[u_{n}=1|s_{n-1}=1\right]=P_{0}\end{array}\right)\end{equation}
\begin{equation}
\left(\begin{array}{cc}
\mathbb{P}\left[s_{n}=0|s_{n-1}=0,u_{n}=0\right]=1 & \mathbb{P}\left[s_{n}=1|s_{n-1}=0,u_{n}=0\right]=0\\
\mathbb{P}\left[s_{n}=0|s_{n-1}=0,u_{n}=1\right]=0 & \mathbb{P}\left[s_{n}=1|s_{n-1}=0,u_{n}=1\right]=1\\
\mathbb{P}\left[s_{n}=0|s_{n-1}=1,u_{n}=0\right]=0 & \mathbb{P}\left[s_{n}=1|s_{n-1}=1,u_{n}=0\right]=1\\
\mathbb{P}\left[s_{n}=0|s_{n-1}=1,u_{n}=1\right]=1 & \mathbb{P}\left[s_{n}=1|s_{n-1}=1,u_{n}=1\right]=0\end{array}\right)\end{equation}
the first matrix is a Bernoulli process transition matrix, the second
shows the dependence between the morphic fundaments ($u_{n}$) and
the previous chain object ($s_{n-1}$) that enact the Bernoulli-type
independence in $S$, the third matrix expresses the deterministic
relation that characterizes the morphic representation of the chain
in $2^{\omega}$.

This constitutes an example of chaos in category computation, such
that one obtains an emergent aleatorial dynamics that is computable,
in the categorial sense, and all the randomness, in the output string
$S$, is removed by the knowledge of the input string and the previous
element $s_{n-1}$.

A statistically independent random-like object sequence, in an infinite
binary morphic chain in $\mathcal{C}_{2}$, may thus be proven to
result from a deterministic dynamical structure, computable through
a categorial automaton in $2^{\omega}$ with a low complexity rule.
Versions of infinite binary morphic chains in $\mathcal{C}_{2}$ are
well explored within dynamical systems theory, in regards to the future
itinerancies of a chaotic orbit whose dynamics is mapped to a binary
symbolic system, which provides for a further piece of conectivity
between chaos, in binary category computation, and continuous state
chaos.

\section{Conclusions}

Category computation theory provides an entry point for a systemic
logic that allows one to connect different fields of mathematics with
effectivenesses towards systems science, dynamical systems theory
and the complexity sciences.

There are two main issues that come interconnected and that provide
for the main points of the present work's argument: the emergence
issue and the algorithmic incompressibility/morphic compressibility
issue. Both, combined, provide for a result that allows one to address
a systemically generated randomness, whose order is accessible and
expressible in small length morphic programs, allowing one to address
algorithmically incompressible strings as being the result of a chaotic
dynamics, generated by a morphism that connects an input (operator)
string and an output string in $2^{\omega}$, the later (string) expressing
the object sequence of an infinite chain in $\mathcal{C}_{2}$.

The emergence issue can be placed as such: for binary categorial automata
in $2^{\omega}$, the fundament of the morphic connection between
an \emph{origin} (\emph{input}) and a \emph{target} (\emph{output})
is an emergent property from local morphic connectivities.

By describing the local connectivities, it becomes possible to address
that fundament in terms of a string that introduces a rule for how
the local morphic connectivities are built. This is the first main
point for the current paper's argument.

The second main point comes from the fact that one can express an
infinite morphic chain in $\mathcal{C}_{2}$, as a single morphism
in $2^{\omega}$, which, in turn, allows one to write a small morphic
program string of twelve bits long, which describes the local connectivities
that generate the infinite morphic chains of $\mathcal{C}_{2}$. This
allows one to show how even algorithmically random binary strings
can be the result of a small length morphic program, such that the
algorithmically random sequence of objects is the result of a systemic
computation which can be expressed by a low length program in $2^{\omega}$.

Even if the dynamics produces a pattern that, considered independently
of the generative systemic fundament, is algorithmically incompressible,
once the generative systemic fundament is addressed, one, thus, recovers
a morphic computational compressibility, expressing a small length
rule that explains how an input generates an output.

Given the fact that there is a one-to-one and onto correspondence
between the infinite binary strings and the object sequences of $\mathcal{C}_{2}$,
it turns out that, treating the infinite binary strings as object
sequences of infinite morphic chains in $\mathcal{C}_{2}$, and then
re-expressing the later as a single morphism in $2^{\omega}$, allows
one to address a categorial generative framework for algorithmically
random infinite binary strings that places these in the realm of binary
categorial chaotic dynamics, generated by small length morphic programs.
As is proved in \emph{section 3.}, all infinite binary strings are
morphically compressible (\emph{theorem 1.}), thus, all algorithmically
random infinite binary strings are morphically compressible.

This main argument, laid out in the present work, was completed by
an exemplification through a probabilistic analysis that addressed
the simultaneous determinism and randomness of an infinite binary
object sequence of an infinite morphic chain in $\mathcal{C}_{2}$,
statistically equivalent to an infinite sequence of independent Bernoulli
random trials.

\end{document}